\documentclass[11pt,reqno]{amsart}

\usepackage{latexsym}
\usepackage{amsfonts}
\usepackage{amsmath}
\usepackage{graphicx}
\usepackage{color}
\usepackage{url}

\usepackage{amssymb,amsmath}

\def\bl{\begin{lemma}}
\def\el{\end{lemma}}
\def\bth{\begin{theorem}}
\def\eth{\end{theorem}}
\def\bc{\begin{corollary}}
\def\ec{\end{corollary}}
\def\bcj{\begin{conjecture}}
\def\ecj{\end{conjecture}}
\def\bpr{\begin{proposition}}
\def\epr{\end{proposition}}
\def\bde{\begin{definition}}
\def\ede{\end{definition}}
\def\E{\mathbb{E}}

\newcommand{\be}{\begin{eqnarray}}
\newcommand{\ee}{\end{eqnarray}}
\newcommand{\eps}{{\mbox{$\epsilon$}}}

\newcommand{\Z}{{\mathbb Z}}

\newcommand{\EE}{{\mathcal E}}

\newcommand{\diam}{{\rm diam}}
\newcommand{\with}{\hbox{ {\rm with} }}
\renewcommand{\and}{\hbox{ {\rm and} }}

\newcommand{\off}{\hbox{ {\rm off} }}
\newcommand{\on}{\hbox{ {\rm only on} }}
\newcommand{\inn}{\hbox{ {\rm in} }}

\newcommand{\Tm}{T_{{\rm mix}}}
\newcommand{\C}{{\mathcal{C}}}
\newcommand{\prob}{\mbox{\bf P}}

\newcommand{\lr}{\leftrightarrow}
\newcommand{\lrlong}{\longleftrightarrow}
\newcommand{\lrr}{\stackrel{\,\, r}{\leftrightarrow}}

\newcommand{\bcr}{B}

\newtheorem{theorem}{Theorem}[section]
\newtheorem{definition}{Definition}[section]
\newtheorem{lemma}[theorem]{Lemma}

\newtheorem{corollary}[theorem]{Corollary}
\newtheorem{proposition}[theorem]{Proposition}
\newtheorem{conjecture}[theorem]{Conjecture}

\theoremstyle{definition}
\numberwithin{equation}{section}

\input epsf.sty

\begin{document}
\title{A note about critical percolation on finite graphs}
\author{Gady Kozma} \author{Asaf Nachmias}

\begin{abstract} In this note we study the geometry of the largest component $\C_1$ of critical percolation on a finite graph $G$ which satisfies the finite triangle condition, defined by Borgs et al. in \cite{BCHSS1}. There it is shown that this component is of size $n^{2/3}$, and here we show that its diameter is $n^{1/3}$ and that the simple random walk takes $n$ steps to mix on it. By \cite{BCHSS2}, our results apply to critical percolation on several high-dimensional finite graphs such as the finite torus $\Z_n^d$ (with $d$ large and $n\to \infty$) and the Hamming cube $\{0,1\}^n$.
\end{abstract}


\maketitle


\section{{\bf  Introduction}}

Given a graph $G=(V,E)$ and $p\in [0,1]$, the probability measure $\prob_p$ on subgraphs of $G$ is obtained by independently retaining each edge with probability $p$ and deleting it with probability $1-p$. Write $G_p$ for the resulting graph and call retained edges {\em open} and deleted edges {\em closed}.
For two vertices $x,y\in V$ and $p \in [0,1]$ the triangle diagram $\bigtriangledown_p(x,y)$ is defined by
$$ \bigtriangledown_p(x,y) = \sum_{u,v \in V} \prob_p(x \lr u) \prob_p(u \lr v) \prob_p(v \lr y) \, ,$$
where $x \lr u$ denotes the event that there exists an open path connecting $x$ to $u$. Given a transitive graph $G$ and $\lambda > 0$, the critical percolation probability $p_c=p_c(G, \lambda)$, defined in \cite{BCHSS1}, is the unique solution to
\be\label{defpc} \E _{p_c} |\C(x)| = \lambda n^{1/3} \, ,\ee
where $x \in G$ is a vertex and $\C(x)$ is the connected component containing $x$ (due to transitivity, the choice of $x$ is arbitrary). The {\em finite triangle condition}, also defined in \cite{BCHSS1}, asserts that
\be \label{triangle} \bigtriangledown_{p_c}(x,y) \leq {\bf 1}_{\{x=y\}} + a_0 \, , \ee for sufficiently small $a_0$ and any $x,y\in G$. We write $\C_1$ for the largest component of $G_p$ and denote by {\rm
diam}$(\C_1)$ and by $T_{{\rm mix}}(\C_1)$ the diameter (maximal
graph distance) of $\C_1$ and the mixing time of the lazy simple
random walk on $\C_1$, respectively (see \cite{NP3} for a
definition). The main theorem of this note is the following.

\begin{theorem} \label{largest} There exists some $a_0 \leq 1/4$ such that for any $\lambda > 0$ and $A < \infty$ the following holds. For any finite transitive graph $G$ of degree $d$ for which (\ref{triangle}) holds with $a_0$, for any $p$ satisfying
$$ p \in [p_c - Ad^{-1}n^{-1/3}, p_c + Ad^{-1}n^{-1/3}] \, ,$$
and any $\eps>0$, there exists $B=B(\eps, \lambda, A)>0$ such that
\begin{enumerate}
\item $\prob_p \big ( \diam(\C_1) \in [B^{-1}n^{1/3}, Bn^{1/3}] \big ) \geq 1-\eps $,
\item $\prob_p \big ( \Tm(\C_1) \in [B^{-1}n, Bn] \big ) \geq 1-\eps $.
\end{enumerate}
\end{theorem}
\noindent The size of $\C_1$ was determined in \cite{BCHSS1}, where the authors prove that under the conditions of Theorem \ref{largest} we have
\be\label{volest}\prob_p \big ( |\C_1| \in [B^{-1}n^{2/3}, Bn^{2/3}] \big ) \geq 1 - \eps \, .\ee

\noindent {\bf Remark.} It is proven in \cite{BCHSS2} that the conditions of Theorem \ref{largest} hold for various high-dimensional graphs such as the torus $\Z_n^d$ when $d$ is large but fixed and $n\to \infty$ or the Hamming hypercube $\{0,1\}^m$. \\

\subsection{Discussion.}
When $G$ is an infinite transitive graph, the {\em triangle condition} states that
\be\label{inftriangle} \bigtriangledown_{p_c}(x,x) < \infty \, ,\ee
where $p_c$ is the critical percolation probability $p_c = \sup\{p:\prob_p(|\C(x)|=\infty)=0\}$. This condition was suggested by Aizenman and Newman \cite{AN} as an indicator for tree-like behavior of critical percolation on infinite graphs. See \cite{AB} and \cite{AN} for further details, and \cite{HaS} for a proof that (\ref{inftriangle}) holds for lattices in high dimensions. When $G$ is a finite transitive graph there is no infinite cluster, and (\ref{inftriangle}) holds for any $p$. The critical phenomenon still occurs, except that the role of the infinite cluster is played by the {\em largest component} $\C_1$ of $G_p$. However, it is not quite clear what is the correct finite case analogue of $p_c$. See \cite{BCHSS1} and \cite{NP3} for further discussion of this topic.

In this note we use the definition of $p_c$ for finite graphs (\ref{defpc}) and assume the finite analogue of the triangle condition (\ref{triangle}), both given by \cite{BCHSS1}. These are expected to represent the actual phase transition {\em only} in the mean-field case, that is, when critical percolation on $G$ behaves as it does on the complete graph $K_n$. Indeed, the estimate (\ref{volest}), proved in \cite{BCHSS1}, asserts that with this $p_c$ the finite triangle condition implies that $\C_1$ is of size $n^{2/3}$ and that this continues for the entire {\em scaling window} $p=p_c(1+O(n^{-1/3}))$. This is precisely what occurs when the underlying graph is $K_n$, see \cite{ER} and \cite{B1}.

In this note we examine further geometric properties of the largest component. In particular, we show that its diameter is of order  $n^{1/3}$ and that simple random walk takes $n$ steps to mix on it. These results are a direct corollary of the following three theorems.
\begin{enumerate}

\item Theorem 1.3 of \cite{BCHSS1}, which guarantees that $|\C_1|$ is of order $n^{2/3}$ when $p$ is in the critical window and the finite triangle condition holds, and
\item Theorem \ref{exp} of this note, asserting that under the finite triangle condition the intrinsic metric critical exponents attain their mean-field values, and
\item Theorem 2.1 of \cite{NP3}, stating that given Theorem \ref{exp}, {\em if} the largest component $\C_1$ is of order $n^{2/3}$, then its diameter is $n^{1/3}$ and its mixing time is $n$.
\end{enumerate}
The proofs in this paper are adaptations of the proofs in \cite{KN}, which relied on (\ref{inftriangle}), to the finite setting. However, a new argument was needed to cover the entire scaling window, that is, when $p\in[p_c, p_c(1+ O(n^{-1/3}))]$.

\subsection{The intrinsic metric critical exponents.}
Let $G$ be a graph and write $G_p$ for the result of $p$-bond
percolation on it. Write $d_{G_p}(x,y)$ for the length of the
shortest path between $x$ and $y$ in $G_p$, or $\infty$ if there
is no such path. We call $d$ the {\em intrinsic metric} on $G_p$.
Define the random sets
\begin{align*}
\bcr_p(x,r;G) &= \{ u : d_{G_{p}}(x,u) \leq r \} \, , \\
\partial \bcr_p(x,r;G,p) &= \{ u : d_{G_{p}}(x,u) = r \} \, .
\end{align*}

Define now the event
\begin{align*}
H_p(r;G) &= \Big \{ \partial B_p(0,r;G) \neq \emptyset \Big \} \, , \\
\intertext{and finally define} \Gamma_p(r; G)&=\sup_{G' \subset G }
\prob (H_p(r;G')) \, .
\end{align*}

The main result of this note is the following.

\begin{theorem}\label{exp} Let $\lambda > 0$ and $A < \infty$. For any finite transitive graph $G$ of degree $d$ for which (\ref{triangle}) hold for some $a_0 \leq 1/4$ and any $p$ satisfying
$$ p \in [p_c - Ad^{-1}n^{-1/3}, p_c + Ad^{-1}n^{-1/3}] \, ,$$
we have
\begin{enumerate}
\item \label{enu:EB0r} \qquad $\E |\bcr_p(0,r; G)| \leq C r$
\item \label{enu:Gam<} \qquad $ \Gamma_p(r; G) \leq {C \over r} \, .$
In particular, $\prob_p(H(r)) \leq {C \over r}$,
\end{enumerate}
where $C>0$ is a constant depending only on $\lambda$ and $A$, but not on $n$ and $d$.
\end{theorem}

Theorem \ref{largest} is an immediate corollary of Theorem \ref{exp} and the two following theorems, which we quote here for the sake of completeness. The first is the aforementioned estimate on the volume of $\C_1$.

\begin{theorem} [Theorem 1.3 of \cite{BCHSS1}] \label{bchss}
There exists some $a_0 \leq 1/4$ such that for any $\lambda > 0$ and $A < \infty$ the following holds. For any finite transitive graph $G$ of degree $d$ for which (\ref{triangle}) holds with $a_0$ and any $p$ satisfying
$$ p \in [p_c - Ad^{-1}n^{-1/3}, p_c + Ad^{-1}n^{-1/3}] \, ,$$
we have that for any $\eps>0$ there exists $B=B(\eps, \lambda, A)>0$ such that
$$\prob_p \big ( |\C_1| \in [B^{-1}n^{2/3}, Bn^{2/3}] \big ) \geq 1 - \eps \, .$$
\end{theorem}

The second is a theorem estimating the diameter and mixing time of the cluster.

\begin{theorem} [Theorem 2.1 of \cite{NP3}] \label{nperes} Let $G=(V,\EE)$ be a graph and $p\in[0,1]$. Suppose that for some constant $C$ and all vertices $v\in V$ we have that
\begin{enumerate}
\item \qquad $\E |\EE(\bcr_p(0,r; G))| \leq C r$
\item \qquad $ \Gamma_p(r; G) \leq {C \over r} \, ,$
\end{enumerate}
where $\EE(\bcr_p(0,r; G))$ denotes the set of open edges with two endpoints in $\bcr_p(0,r; G)$. Then:
\begin{enumerate}
 \item $\prob_p \big ( \exists \C \with |\C|\geq \beta n^{2/3} \and \diam(\C) \not \in [B^{-1} n^{1/3}, Bn^{1/3}] \big ) \leq O(B^{-1})$ ,
 \item $\prob_p \big ( \exists \C \with |\C|\geq \beta n^{2/3} \and \Tm(\C) \geq Bn \big ) \leq O(B^{-1/2})$ ,
 \item $\prob_p \big ( \exists \C \with |\C|\geq \beta n^{2/3} \and \Tm(\C) \leq B^{-1}n \big ) \leq O(B^{-1/13})$ ,
\end{enumerate}
where the constants implicit in the $O$-notation depend only on $C$ and $\beta$.
\end{theorem}

\noindent {\bf Proof of Theorem \ref{largest}}. Indeed, the theorem follows directly from Theorems \ref{exp} \ref{bchss} and \ref{nperes}, except that Theorem \ref{exp} gives that $\E |\bcr_p(0,r; G)| \leq C r$ and we need to verify the same for $\E |\EE(\bcr_p(0,r; G))|$ in order to use Theorem \ref{nperes}. Indeed, consider ``exploring'' the levels $\partial \bcr_p(0,k; G)$ level by level for $k=1, \ldots, r$. At the end we discovered a spanning tree on the vertices of $\bcr_p(0,r; G)$ and since the degree is $d$, the number of extra edges in this ball can be bounded above by a random variable $Z$ distributed as Bin$(d|\bcr_p(0,r; G)|,p)$. Thus, if we condition on $|\bcr_p(0,r; G)|$, then the number of edges $|\EE(\bcr_p(0,r; G))|$ can be stochastically bounded above by $|\bcr_p(0,r; G)| + Z$. We now appeal to Theorem 1.1 of \cite{BCHSS1} which implies that in our setting $dp \leq C$ for some constant $C=C(\lambda, A)>0$.
\qed

\subsection{Related work.} Consider the finite torus $\Z_n^d$ when $d$ is large but fixed and $n\to \infty$ and perform percolation with $p=p_c(\Z^d)$, that is, the critical probability for the infinite lattice $\Z^d$. Van der Hofstad and Heydenreich \cite{HH1, HH2} show that the largest component arising is of size $n^{2/3}$ and that $|p_c(\Z^d) - p_c(\Z_n^d, \lambda)|=O(n^{-1/3})$. The authors then appeal to the infinite version of Theorem \ref{exp} in \cite{KN} and to Theorem \ref{nperes} to obtain the estimates for the diameter and mixing time of the largest component when $p=p_c(\Z^d)$.

In fact, some of the estimates of this paper exist in \cite{HH2} in the same level of generality. In \cite{HH2} the authors adapt the proofs of \cite{KN}, as we do here, to the finite case and establish Theorem \ref{exp}, albeit only for $p \leq p_c(G,\lambda)$. Our main effort here is handling the entire scaling window $p = (1+\Theta(n^{-1/3}))p_c$  which includes the case $p>p_c$. The latter is not covered in the results of \cite{HH1} and \cite{HH2} and requires a new {\em sprinkling} argument, see the proof of part (\ref{enu:EB0r}) of Theorem \ref{exp}. Both in \cite{HH2} and here the diameter and mixing time estimates are direct corollaries.

\section{Proofs}

From now on, we assume the condition of Theorem \ref{exp}, that is, we assume that for a given $G$ and $\lambda>0$ (which determine $p_c$ by (\ref{defpc})) the triangle condition (\ref{triangle}) holds with $a_0 \leq 1/4$. Denote
$$ G(r) = \E |\bcr_{p_c}(0,r; G)| \, .$$
The main part of the proof is the following lemma.
\begin{lemma} \label{volrecur} Under the setting of Theorem \ref{exp} with $p=p_c$ we have that for all $r$
$$ G(2r) \geq {G(r)^2 \over 4r} \, .$$
\end{lemma}

\noindent Let us first see how to use the lemma \ref{volrecur}.

\noindent {\bf Proof of part (\ref{enu:EB0r}) of Theorem
\ref{exp}.} We first prove the assertion for $p\leq p_c$ and then use this estimate together with a sprinkling argument to prove the assertion for $p \in [p_c, p_c(1+O(n^{-1/3})]$.

\noindent{\em Proof for $p \leq p_c$.} Note that the random variable $|\bcr_{p}(0,r; G)|$ is monotone with respect to adding edges, hence it suffices to prove the statement for $p=p_c$. We prove that $G(r) \leq 8 r$ for all $r$. Assume by contradiction that there exists $r_0$ such that $G(r_0) \geq 8 r_0$. Under this assumption, we prove by induction that for any
integer $k \geq 0$ we have $G(2^k r_0) \geq 8^{k+1}r_0$. The
case $k=0$ is our assumption and for $k \geq 1$ Lemma \ref{volrecur} gives that
$$ G(2^{k+1} r_0) \geq {G (2^k r_0)^2 \over 4\cdot  2^k r_0} \geq 8^{k+2} r_0 \,
,$$ where in the last inequality we used the induction hypothesis
and the fact that $8^k \geq 4 \cdot 2^k$ for any $k \geq 1$. This completes our induction.

We have now arrived at a contradiction, since on the right hand side we have a sequence going to $\infty$ with $k$, and on the left hand side we have a sequence in $k$ which is bounded by $|G|$. We learn that $\E |\bcr_{p_c}(0,r; G)| \leq 8r$ for all $r$. This concludes the case $p \leq p_c$.

\noindent {\em Proof for $p \in [p_c, p_c + Ad^{-1}n^{-1/3}]$.} Again, due to monotonicity it suffices to prove for $p=p_c+\delta$ where $\delta=Ad^{-1}n^{-1/3}$. We recall the usual simultaneous coupling of $\prob_{p_c}$ and $\prob_{p_c+\delta}$. In this coupling we assign to each edge $e$ an i.i.d uniform $[0,1]$ random variable $X_e$ and an edge $e$ is declared $p$-open if $X_e \leq p$, for any $p\in[0,1]$. In this coupling, we say that an edge $e$ is  {\em sprinkled} if $X_e \in [p_c, p_c+\delta]$. For integers $r,m$ we write $0 \stackrel{r,m}{\lrlong} x$ for the event that $0$ is connected to $x$ in a simple path of length at most $r$ which have all its edges $p_c$-open except for precisely $m$ edges, which are sprinkled. It is clear that \be\label{whostheboomking?} \bcr_{p_c+\delta}(0,r; G) \subset \bigcup_{m \geq 0} \{x : 0 \stackrel{r,m}{\lrlong} x\} \, .\ee
We now prove by induction on $m$ that
$$ \E \big | \{x : 0 \stackrel{r,m}{\lrlong} x\}\big | \leq 8r \cdot (16 d \delta r)^m \, .$$
For $m=0$ we have that $\E| \{x : 0 \stackrel{r,0}{\lrlong} x\} | = \E |\bcr_{p_c}(0,r; G)| \leq 8r$. Now fix $m \geq 1$. If $0 \stackrel{r,m}{\lrlong} x$ occurs, then there exist an edge $(y,y')$ in $G$ such that the event
$$ \{0 \stackrel{r,m-1}{\lrlong} y\} \circ \{(y,y') \textrm{ is sprinkled}\} \circ \{y' \stackrel{r,0}{\lrlong} x\} \, ,$$
occurs. Summing over $(y,y')$ and applying the BK inequality gives that
$$ \E \big | \{x : 0 \stackrel{r,m}{\lrlong} x\}\big | \leq \E \big | \{y : 0 \stackrel{r,m-1}{\lrlong} y\}\big |\cdot 2d\delta \cdot \E |\bcr_{p_c}(0,r; G)| \, ,$$ and the induction assertion follows.  Putting this into (\ref{whostheboomking?}) and summing over $m$ shows that for $r \leq {n^{1/3} \over 32 A}$ we have
$$ \E \big | \bcr_{p_c+\delta}(0,r; G) \big | \leq 16r \, .$$
For $r \geq {n^{1/3} \over 32 A}$ we simply bound $|\bcr_{p}(0,r; G)|\leq |\C(0)|$, and part (iii) of Theorem 1.3 of \cite{BCHSS1} gives that $\E_p |\C(0)| \leq C(\lambda, A) n^{1/3}$, concluding the proof. \qed \\

\noindent We now turn to the proof of Lemma \ref{volrecur}.

\begin{lemma} \label{revBK} Let $0\in G$ be a vertex. Under the conditions of Theorem \ref{exp} with $p=p_c$ we have
$$ \sum _{x,y \in G} \prob_p \big ( \{0 \lrr x\} \circ \{x \lrr y\}
\big ) \geq {3 G(r)^2 \over 4} \, .$$
\end{lemma}
\noindent {\bf Proof.} Since $G$ is a finite transitive graph we have that
$$ \sum _{x,y \in G} \prob_p \big ( \{0 \lrr x\} \circ \{x \lrr y\}
\big ) = \sum _{x,y \in G} \prob_p \big ( \{0 \lrr x\} \circ \{0 \lrr y\}
\big ) \, .$$
In fact, this equality holds for any transitive {\em unimodular} graph.
Thus, it suffices to prove that
\be\label{revBK.step0} \sum _{x,y \in G} \prob_p \big ( \{0 \lrr x\} \circ \{0 \lrr y\}
\big ) \geq {3 G(r)^2 \over 4} \, .\ee
For a vertex $z \in G$ we write $\C^{0}(z)$ for
$$ \C^0(z) = \{ u \in G : z \lr u \off 0 \} \, ,$$
where the event $z \lr u \off 0$ means there exists an open path connecting $z$ to $u$ which avoids $0$. We have
$$ \prob_p \big (  \{0 \lrr x\} \circ \{0 \lrr y\} \big ) \geq \prob_p
\Big ( 0 \lrr x \and 0 \lrr y \and \C^0(y) \neq \C^0(x) \Big )
\, .$$ By conditioning on $\C^0(x)$
we get that the right hand side equals
\begin{eqnarray*}
  \sum_{\substack{A\subset G\\
                  0\lrr x \inn A\cup \{0\}, y \not \in A}}
\prob_p(\C^0(x)=A)\prob_p \Big (0 \lrr y \and \C^0(y) \cap A = \emptyset \, \mid \,  \C^0(x) = A \Big ) \, .
\end{eqnarray*}
For $A$ such that $y \not \in A$ we have that
$$\prob_p \big ( \{0 \lrr y\} \and \C^0(y) \cap A = \emptyset \, \mid \,  \C^0(x) = A \big ) = \prob_p \big ( 0 \lrr y \off A \big )$$
where the event $\{0 \lrr y$ off
$A\}$ means that there exists an open path of length at most
$r$ connecting $0$ to $y$ which avoids the vertices of $A$. At
this point we can remove the condition that $y \not \in A$
since in this case the event $\{0 \lrr y$ off $A\}$ is
empty. We get
\begin{equation}\label{condcu}
\prob_p \big (  \{0 \lrr x\} \circ \{0 \lrr y\} \big ) \geq
  \sum_{\substack{A\subset G\\
                  0\lrr x \inn A\cup \{0\}}}
\prob_p(\C^0(x)=A) \prob_p \big ( 0 \lrr y \off A \big ) \, .
\end{equation}
Now, since
$$
\prob_p ( 0 \lrr x) \prob_p (0 \lrr y) =
\sum_{\substack{A \subset G\\
                0\lrr x \inn A\cup \{0\}}}
\prob_p(\C^0(x)=A)  \prob_p (0 \lrr y) \, ,$$
we deduce by (\ref{condcu}) that \begin{multline}
\label{revBKstp1} \prob_p \big (  \{0 \lrr x\} \circ \{0 \lrr y\} \big ) \geq \prob_p ( 0 \lrr x) \prob_p (0 \lrr y) \\
- \sum _{\substack{A \subset G\\
                   0\lrr x \inn A\cup \{0\}}}
\prob_p(\C^0(x)=A)) \prob_p \big ( 0 \lrr y \on A \big ) \,
,\end{multline} where the event $\{0 \lrr y \on A\}$ means that
there exists an open path between $0$ and $y$ of length at most
$r$ and any such path must have a vertex in $A$. If $\{0 \lrr y \on A\}$ occurs, then there must exists $z \in A$ such that $\{0 \lr z\} \circ \{z \lrr y\}$ occurs. Hence for any $A
\subset G$ we have
$$ \prob_p \big ( 0 \lrr y \on A \big ) \leq \sum _{z \in A}
\prob_p \big ( \{0 \lr z\} \circ \{z \lrr y\} \big ) \, .$$
Putting this into the second term of the right hand side of
(\ref{revBKstp1}) and changing the order of summation gives that
we can bound this term from above by
$$\sum_{z \in G\setminus \{0\}} \prob_p\big(0 \lrr x \, , \, 0 \lr
z \big ) \prob_p \big ( \{0 \lr z\} \circ \{z \lrr y\} \big ) \,
,$$ where the sum is over $z \neq 0$ since $0 \not \in A$.
If $0 \lrr x$ and $0 \lr z$, then there exists $z'$ such
that the event $\{0 \lr z'\} \circ \{z' \lr z\} \circ \{z' \lrr x\}$ occurs. Using the BK inequality we bound this sum
above by
$$ \sum_{z \in G \setminus \{0\}, z'\in G} \prob_p(0 \lr z') \prob_p(z' \lr z) \prob_p (z'
\lrr x) \prob_p(0 \lr z) \prob_p(z \lrr y) \, .$$
We sum this over $x$ and $y$ and use (\ref{revBKstp1}) to get that
$$ \sum _{x,y\in G} \prob_p \big (  \{0
\lrr x\} \circ \{0 \lrr y\} \big ) \geq G(r)^2 -
G(r)^2 \sum _{z \in G \setminus \{0\}, z'\in G} \prob_p(0 \lr z') \prob_p(z' \lr z) \prob_p(z \lr 0) \, .$$
The finite triangle condition (\ref{triangle}) and the fact that $z \neq 0$ (excluding from the sum the term $z=z'=0$, which equals $1$) implies that
$$ \sum _{z \in G \setminus \{0\}, z'\in G} \prob_p(0 \lr z') \prob_p(z' \lr z) \prob_p(z \lr 0) \leq a_0 \, ,$$
and since $a_0 \leq 1/4$ we get the assertion of the lemma. \qed

\begin{figure}
\begin{picture}(0,0)%
\includegraphics{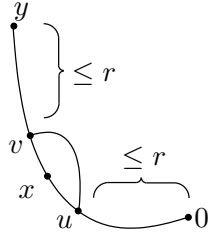}%
\end{picture}%
\setlength{\unitlength}{4144sp}%
\begingroup\makeatletter\ifx\SetFigFont\undefined%
\gdef\SetFigFont#1#2#3#4#5{%
  \reset@font\fontsize{#1}{#2pt}%
  \fontfamily{#3}\fontseries{#4}\fontshape{#5}%
  \selectfont}%
\fi\endgroup%
\begin{picture}(1251,1435)(4395,-3521)
\put(4456,-3301){\makebox(0,0)[lb]{\smash{{\SetFigFont{11}{13.2}{\rmdefault}{\mddefault}{\updefault}{\color[rgb]{0,0,0}$x$}%
}}}}
\put(4681,-3481){\makebox(0,0)[lb]{\smash{{\SetFigFont{11}{13.2}{\rmdefault}{\mddefault}{\updefault}{\color[rgb]{0,0,0}$u$}%
}}}}
\put(5508,-3476){\makebox(0,0)[lb]{\smash{{\SetFigFont{11}{13.2}{\rmdefault}{\mddefault}{\updefault}{\color[rgb]{0,0,0}$0$}%
}}}}
\put(4428,-2206){\makebox(0,0)[lb]{\smash{{\SetFigFont{11}{13.2}{\rmdefault}{\mddefault}{\updefault}{\color[rgb]{0,0,0}$y$}%
}}}}
\put(4395,-3029){\makebox(0,0)[lb]{\smash{{\SetFigFont{11}{13.2}{\rmdefault}{\mddefault}{\updefault}{\color[rgb]{0,0,0}$v$}%
}}}}
\put(5073,-3103){\makebox(0,0)[lb]{\smash{{\SetFigFont{11}{13.2}{\rmdefault}{\mddefault}{\updefault}{\color[rgb]{0,0,0}$\leq r$}%
}}}}
\put(4775,-2599){\makebox(0,0)[lb]{\smash{{\SetFigFont{11}{13.2}{\rmdefault}{\mddefault}{\updefault}{\color[rgb]{0,0,0}$\leq r$}%
}}}}
\end{picture}%
\caption{\label{cap:overcounted}The couple $(x,y)$ is over-counted.}
\end{figure}
\medskip

\noindent {\bf Proof of Lemma \ref{volrecur}.} We start with a
definition. We say two vertices $x,y \in
G$ are {\em over-counted} if there exists $u,v \in G\setminus \{x\}$ with
such that
$$ \{0 \lrr u\} \circ \{v \lr x\} \circ \{x \lr u\} \circ \{u
\lr v\} \circ \{v \lrr y \} \, ,$$ see figure \ref{cap:overcounted}.
Denote by
$N$ the quantity
$$ N = \Big | \Big \{ (x,y) \, : \, \{0 \lrr x\} \circ \{x \lrr
y\} \and (x,y) \hbox{ {\rm are not over-counted}} \Big \} \Big
| \, .$$
We claim that
\be \label{determ} N \leq 2r |\bcr(0,2r)| \, .\ee Indeed, this deterministic claim follows by
observing that if $y \in \bcr(0,2r)$ and $\gamma$ is an open
simple path of length at most $2r$ connecting $0$ to $y$, then
for any $x \in G\setminus \gamma$
satisfying $\{0 \lrr x\}\circ \{x \lrr y\}$ the pair $(x,y)$ is
over-counted. To see this, let $\gamma_1$ and $\gamma_2$ be
disjoint open simple paths of length at most $r$ connecting $0$ to
$x$ and $x$ to $y$ respectively and take $u$ to be the last point
on $\gamma \cap \gamma_1$ and $v$ the first point on $\gamma \cap
\gamma_2$ where the ordering is induced by $\gamma_1$ and $\gamma_2$
respectively. Hence the map
$(x,y)\mapsto y$ from $N$ into $B(0,2r)$ is at most $2r$ to $1$,
which shows (\ref{determ}).

We now estimate $\E N$. For any $(x,y)$ the BK inequality implies that the probability that $(x,y)$ is
over-counted is at most
$$ \sum_{\substack{u \neq x,v \neq x}} \prob_p(0 \lrr u) \prob_p(v \lrr y) \prob_p(x \lr v) \prob_p(v \lr u) \prob_p(u \lr x) \, .$$
Denote by $\varphi_x$ a graph automorphism taking $x$ to $0$ (which exists by transitivity). Since $G$ is transitive, the last sum equals to
$$
 \sum_{\substack{u' \neq 0, v' \neq 0}}
  \prob_p(\varphi_x(0) \lrr u') \prob_p(v' \lrr \varphi_x(y)) \prob_p (0 \lr v') \prob_p(v' \lr u') \prob_p (u' \lr 0) \, .
$$
We sum this over $y$ and then over $x$
to get that \begin{align*} \E \Big [ \big | \big \{ (x,y)
\hbox{ {\rm are over-counted}}\big \} \big | \Big ] &\leq
G(r)^2
\sum_{\substack{u' \neq 0, v'\neq 0}} \prob_p (0 \lr v') \prob_p(v' \lr u') \prob_p (u' \lr 0) \, .
\end{align*}
The triangle condition (\ref{triangle}) implies that the sum on the right hand side is at most ${1 \over 4}$. Together with Lemma \ref{revBK} this implies that $\E N \geq {G(r)^2/2}$, which concludes the proof of the lemma by (\ref{determ}). This also concludes the proof of part (i) of Theorem \ref{exp}.
\qed \\

We proceed with the proof of part (ii) of Theorem \ref{exp}. The following is a more general theorem, whose proof follows verbatim the proof of part (ii) of Theorem 1.2 of \cite{KN}. We state it here and omit the proof.

\begin{theorem} \label{intonearm} Let $G$ be a transitive graph (finite or infinite) and $0$ is an arbitrary vertex and let $p\in [0,1]$ be given. Assume there exists $C_1>0$ such that
\be\label{voldec} \prob_p( |\C(0)| \geq k ) \leq {C_1 \over \sqrt{k}} \, ,\ee
for all $k \geq 1$. Then there exists a constant $C_2 = C_2(C_1)>0$ such that
$$ \Gamma_p(r) \leq {C_2 \over r} \, ,$$
for all $r \geq 1$.
\end{theorem}

\noindent {\bf Proof of part (\ref{enu:Gam<}) of Theorem
\ref{exp}.} This follows from Theorem \ref{intonearm}.  The estimate (1.19) of \cite{BCHSS1} shows (\ref{voldec}) for $k=O(n^{2/3})$ and for larger $k$'s, this follows from the estimate (1.22) of \cite{BCHSS1} and Proposition $5.1$ of \cite{AN}. \qed \\

\noindent {\bf Acknowledgements:} We wish to thank Markus Heydenreich and Remco van der Hofstad for useful comments and for suggesting the use of Proposition $5.1$ of \cite{AN} in the proof of part (\ref{enu:Gam<}) of Theorem \ref{exp}.

\vspace{.05 in}\noindent
{\bf Gady Kozma}: \texttt{gady.kozma(at)weizmann.ac.il} \\
The Weizmann Institute of Science, \\
Rehovot POB 76100, \\
Israel.

\medskip \noindent
{\bf Asaf Nachmias}: \texttt{asafn(at)microsoft.com} \\
Microsoft Research,
One Microsoft way,\\
Redmond, WA 98052-6399, USA.

\end{document}